\documentclass[a4paper,11pt]{article}
\usepackage[french]{babel}
\usepackage{amssymb}

\def \fu #1#2{{\bf\uppercase{#1\footnotesize{#2}}}}
\def \th #1#2#3{\noi\fu{#1}{#2} {\bf #3.} ---}
\def \thf #1#2{\noi\textsc{#1 #2.} ---}

\def \p {{\mathbb P}}

\def \fl {\longrightarrow}
\def \ot {\otimes}

\def \vs {\vskip}

\def \oo {{\cal O}}
\def \dm {{\textit{D\'emonstration}}}

\def \noi {\noindent}
\def \U {{\bf U}}

\def \E {{\cal E}}

\def \op {\oo_{X}}

\def \h {\underline{\bf Hom}}

\def \h {\underline{\bf Hom}_{\op}}

\def \M {{\bf M}}
\def \U {{\bf U}}

\def \Q {\mathbb{Q}}
\def \CU {\overline{\U}}
\def \cpi {\overline{\pi}}
\def \bua {\partial\U_1}
\def \bub {\partial\U_2}

\begin{document}

\centerline{\large{\uppercase{\bf{Limites de fibr\'es vectoriels}}}}
\centerline{\large{\uppercase{\bf{dans $\M_{\Q_3}(0,2,0)$}}}}

\vs 0.3 cm

\centerline{\large{Nicolas PERRIN}}


%
%

\vs 1.3 cm

\centerline{\large{\bf Introduction}} 

\vs 0.5 cm

Soit $\Q_3$ la quadrique lisse de dimension 3 de $\p^4$. Dans \cite{OS},
G. Ottaviani et M. Szurek d\'ecrivent l'ouvert $\U$ de
l'espace des modules $\M_{\Q_3}(0,2,0)$ form\'e des fibr\'es vectoriels de
rang 2 et de classes de Chern $(0,2,0)$. Dans cette note, nous d\'ecrivons
l'adh\'erence $\CU$ de $\U$ dans l'espace
$\M_{\p^3}(0,2,0)$. Nous montrons que le bord est form\'e de deux
composantes irr\'eductibles et qu'elles v\'erifient les ``conditions au
bord'' \'enonc\'ees dans \cite{PE}. Plus pr\'ecis\'ement, nous montrons
le 


\vs 0.4 cm

\th{T}{h\'eor\`eme}{1} \textit{La vari\'et\'e $\CU$ est isomorphe \`a
  $\p_{\p^4}(\Lambda^2{\check \Omega^1_{\p^4}})$ (de projection $\pi$ vers
  $\p^4$).
Soient $\bua$ le ferm\'e
  $\pi^{-1}(\Q_3)$ et $\bub\subset\p_{\p^4}(\Lambda^2{\check
  \Omega^1_{\p^4}})$ le fibr\'e en quadriques naturel au
  dessus de $\p^4$. Les ferm\'es $\bua$ et $\bub$ forment les composantes
  irr\'eductibles du bord de $\U$.}

\textit{Si $E$ est un faisceau g\'en\'eral de $\bua$, alors on a la suite
  exacte :
$$0\fl E\fl\oo_{\Q_3}^2\fl\oo_C(1)\fl 0$$
o\`u $C$ est une conique lisse de $\Q_3$. Si $E$ est un faisceau
g\'en\'eral de $\bub$, alors on a la suite  exacte :
$$0\fl E\fl E''\fl \oo_P\fl 0$$
o\`u $P$ est un point de $\Q_3$ et $E''$est r\'eflexif.}


\vs 0.4 cm


Notons $V$ l'espace vectoriel $H^0\oo_{\p^4}(1)$. Si $P$ est un point de
$\p^4$ nous noterons $\p^4\setminus\{ P\}\stackrel{p}{\fl}\p^3$ la
projection et $\p^4\setminus\{ P\}\stackrel{i}{\fl}\p^4$ l'immersion ouverte. 
Notons $K$ le faisceau $\Omega_{\p^4}(1)\boxtimes\oo_{\p^4}(-1)$ sur
$\p^4\times\p^4$. On a la r\'esolution suivante de la diagonale $\Delta$
(voir \cite{OSS}) : 
$$0\fl\Lambda^4{\check K}\fl \Lambda^3{\check K}\fl
\Lambda^2{\check K}\fl {\check
  K}\fl\oo_{\p^4\times\p^4}\fl\oo_{\Delta}\fl 0$$ 
On peut d\'efinir le fibr\'e ${\cal G}$ sur $\p^4\times\p^4$ (muni de $p_1$
et $p_2$) par la suite exacte : 
$$0\fl{\cal G}\fl\Omega^1_{\p^4}(1)\boxtimes\oo_{\p^4}\fl
\oo_{\p^4}\boxtimes\oo_{\p^4}(1)\fl\oo_{\Delta}(1)\fl 0$$ 
Au dessus du point $P\in\p^4$, le faisceau ${\cal
  G}\vert_{p_1^{-1}(P)}$ s'identifie au faisceau $i_*p^*\Omega^1_{\p^3}(1)$.
Le faisceau ${p_1}_*(\h(\oo_{\p^4}\boxtimes\oo_{\p^4}(-1),{\cal G}))$
s'identifie \`a $\Lambda^2(\Omega^1_{\p^4}(1))$. 


Notons alors $X=\p_{\p^4}(\Lambda^2{\check \Omega^1_{\p^4}})$. 
Sur $X\times\p_4$ (nous notons encore $p_1$ et $p_2$ les
projections), on peut d\'efinir le faisceau sans torsion $\E$ suivant :
$$0\fl\oo_X\boxtimes\oo_{\p^4}(-1)\fl{\cal G}\fl\E\fl 0$$


Nous noterons $(P,s)$ les points de $X$ o\`u $P\in\p^4$ et
$s\in H^0{\cal G}(1)\vert_{p_1^{-1}(P)}$.

\vs 0.4 cm

\th{F}{ait}{2} \textit{On a un isomorphisme $\E\vert_{p_1^{-1}(P,s)}=i_*p^*N$
  o\`u $N\in\M_{\p^3}(0,1,0)$.}

\vs 0.4 cm

Le faisceau $N$ est soit localement libre (instanton de degr\'e 1) soit
donn\'e par le noyau de la fl\`eche de droite dans la suite exacte 
$$0\fl\oo_{\p^3}(-1)\fl\Omega_{\p^3}^1(1)\fl\oo_{\p^3}^2\fl\oo_{L_s}(1)\fl 0$$
o\`u $L_s$ est une droite de $\p^3$. 
Dans ce cas on a 
$$0\fl\E\vert_{p_1^{-1}(P,s)}\fl\oo_{\p^4}^2\fl\oo_{H_s}(1)\fl\oo_P\fl 0$$
o\`u $H_s$ est le plan contenant $P$ se projettant sur $L_s$. Les classes
de Chern de ces faisceaux sont constantes, la famille $\E$ est donc
plate au dessus de $X$. Le faisceau $\E\ot\oo_{X\times\Q_3}$ est encore
sans torsion et plat au dessus de $X$. Ses classes de Chern sont $(0,2,0)$. Il
d\'efinit donc un morphisme 
$$f:X\fl\M_{\Q_3}(0,2,0)$$
qui est birationnel sur $\CU$ ($X\setminus(X_1\cup X_2)$ s'envoie
bijectivement sur $\U$, cf. \cite{OS}).
Nous allons maintenant calculer la restriction de $\E\vert_{p_1^{-1}(P,s)}$
\`a la quadrique $\Q_3$ pour les points de $X_1\cup X_2$. 

\vs 0.4 cm

%
%

\th{P}{roposition}{3} \textit{Le morphisme $f$ est un isomorphisme sur
  $\CU$.}

\vs 0.2 cm

\dm :
Il suffit de v\'erifier qu'il est injectif sur $X_1\cup X_2$, c'est \`a
dire que le faisceau permet de retrouver le point $(P,s)$. Nous noterons
$E$ le faisceau $(\E\vert_{X\times\Q_3})\vert_{p_1^{-1}(P,s)}$

Si $(P,s)\in X_1\setminus (X_1\cap X_2)$, alors on a la suite exacte :
$$0\fl E\fl\oo_{\Q_3}^2\fl
\oo_{H_s\cap\Q_3}(1)\fl 0.$$
La courbe $H_s\cap\Q_3$ est une conique, c'est le lieu singulier du
faisceau. Elle permet de retrouver le plan $H_s$. Les deux sections de
$\oo_{H_s\cap\Q_3}(1)$ d\'efinissent le point $P$ (c'est le lieu d'annulation
de $\oo_{H_s}^2\fl\oo_{H_s}(1)$). La projection de $H_s$ par $P$ d\'efinit
alors la droite $L_s$ qui permet de retrouver la section $s$.

Si $(P,s)\in X_2$, alors la restriction \`a $p_1^{-1}(P,s)\cap\Q_3$ de la
suite exacte de d\'efinition de ${\cal G}$ donne la suite exacte 
$$0\fl{\cal G}\vert_{p_1^{-1}(P,s)\cap\Q_3}\fl H^0\oo_{\p^3}(1)\ot
\oo_{\Q_3}\fl{\cal I}_{P,\p^4}\ot\oo_{\Q_3}\fl 0.$$ 
Cependant le faisceau ${\cal I}_{P,\p^4}\ot\oo_{\Q_3}$ est une extension
de ${\cal I}_{P,\Q_3}$ par $\oo_P$. Ainsi on a la suite exacte :
$$0\fl{\cal G}\vert_{p_1^{-1}(P,s)\cap\Q_3}\fl{\cal G}
\vert_{p_1^{-1}(P,s)\cap\Q_3}''\fl\oo_P\fl 0$$
o\`u le faisceau ${\cal G}\vert_{p_1^{-1}(P,s)\cap\Q_3}''$ est
r\'eflexif. Nous en d\'eduisons la suite exacte 
$$o\fl E\fl E''\fl\oo_P\fl 0$$
o\`u $E''$ est r\'eflexif singulier au point $P$. Le point $P$ est le point
singulier de $E$. Le faisceau $i^*E$ est un fibr\'e vectoriel sur
$\Q_3\setminus\{ P\}$. Le faisceau $p_*i^*E$ est le faisceau
$N\in\M_{\p^3}(0,1,0)$ (tel que $\E\vert_{p_1^{-1}(P,s)}=i_*p^*N$), il
d\'etermine la section $s$. 

Remarquons que si $(P,s)\in X_1\cap X_2$, alors $E$ est donn\'e par la
suite exacte
$$o\fl E\fl\oo_{\Q_3}\oplus{\cal I}_{C,\Q_3}\fl\oo_P\fl 0$$
o\`u $C\subset\Q_3$ est une conique contenant le point $P$.

\vs 0.4 cm

\thf{Remarque}{4} Si on note $\p^4_{\p^4}=\p^4\times\p^4$ et
$\p^3_{\p^4}=\p_{\p^4}(\Omega^1_{\p^4})$, on a alors une
projection et une immersion universelles 
$$p:\p^4_{\p^4}\setminus\Delta\fl\p^3_{\p^4}\ \ \ {\rm et}\ \ \
i:\p^4_{\p^4}\setminus\Delta\fl\p^4_{\p^4}.$$
Notons $j$ l'immersion de $\Q_3$ dans $\p^4$. 
Nous avons montr\'e que le morphisme :
$$\M_{\p^3_{\p^4}}(0,1,0)\stackrel{i_*p^*}{\fl}\M_{\p^4}(0,1,0)
\stackrel{j^*}{\fl}\M_{\Q_3}(0,2,0)$$ 
est un isomorphisme.

\vs 0.4 cm

\thf{Remarque}{5} Il est ici facile de construire des d\'eformations
explicites, on peut alors montrer qu'un \'el\'ement g\'en\'eral du bord est
limite ``r\'eduite'' de fibr\'es vectoriels (au sens de \cite{PE}). 
Les limites ci-dessus v\'erifient les r\'esultats de \cite{PE} : 

Sur $X_1\setminus (X_1\cap X_2)$, le faisceau est localement libre en
  dehors d'un lieu de dimension pure \'egale \`a 1. Il v\'erifie la
  condition du th\'eor\`eme 0.1 : le faisceau
  $\oo_C(1)\ot\oo_C(-\frac{3}{2})$ est une th\'eta-caract\'eristique (le
  faisceau $\oo_C(-\frac{3}{2})$ est une racine de $\omega_{\Q_3}\vert_C$).
Les conditions des th\'eor\`emes 0.3
  et 0.4 sont vides.

Sur $X_2\setminus (X_1\cap X_2)$, le faisceau est singulier un point. De
plus on a la compos\'ee
$${\rm Hom}(E'',\oo_P)\fl{\rm Hom}({\cal
  G}\vert_{p_1^{-1}(P,s)\cap\Q_3}'',\oo_P)\fl{\rm Ext}^1({\cal
  I}_{P,\Q_3},\oo_P)={\rm Ext}^2(\oo_P,\oo_P)$$
et l'image de la surjection de $E''$ dans $\oo_P$ (d\'efinissant $E$) est
l'\'el\'ement de ${\rm Ext}^1({\cal I}_{P,\Q_3},\oo_P)$ d\'efini par le
faisceau ${\cal I}_{P,\p^4}\ot\oo_{Q_3}$. Cet \'el\'ement est nul, c'est
la condition du th\'eor\`eme 0.3.


\vs 0.2 cm

\noi
\textsc{Mathematisches Institut der Universit\"at zu k\"oln}
\vs -0.1 cm
\noi
Weyertal 86-90
\vs -0.1 cm
\noi
D-50931 K\"oln 
\vs -0.1 cm
\noi
email : \texttt{nperrin@mi.uni-koeln.de}


\begin{thebibliography}{99}
\bibitem[OS]{OS} \textit{Giorgio Ottaviani et Micha\l\ Szurek} : On moduli of
  stable $2$-bundles with small Chern classes on $Q\sb 3$. With an appendix
  by Nicolae Manolache. Ann. Mat. Pura Appl. (4) 167 (1994).
\bibitem[OSS]{OSS} \textit{Christian Okonek, Michael Schneider et Heinz
    Spindler} : Vector bundles on complex projective spaces. Progress in
  Mathematics, 3. Birkh\"auser, Boston, Mass. (1980).
\bibitem[P]{PE} \textit{Nicolas Perrin} : D\'eformations de fibr\'es
  vectoriels sur les vari\'et\'es de dimension 3, pr\'eprint.
\end{thebibliography}
\end{document}